\documentclass[reqno, 12pt]{amsart}
\usepackage{amsfonts}
\usepackage{amssymb,amsmath}
\usepackage{amsthm}
\usepackage{upref}
\usepackage{enumerate}
\usepackage{bbm}

\usepackage{amscd}

\makeatletter
\def\LaTeX{\leavevmode L\raise.42ex
    \hbox{\kern-.3em\size{\sf@size}{0pt}\selectfont A}\kern-.15em\TeX}
 \makeatother

\numberwithin{equation}{section}

\newtheorem{lemma}{Lemma}[section]
\newtheorem{theorem}[lemma]{Theorem} 
\newtheorem{corollary}[lemma]{Corollary}

\theoremstyle{definition}

\newcommand{\dist}{\operatorname{dist}}

 \newcommand{\supp}{\operatorname{supp}}
  \newcommand{\e}{\eqref}

\newcommand{\q}{\quad}

\newcommand{\wt}{\widetilde}

\newcommand{\la}{\langle}
\newcommand{\ra}{\rangle}
\newcommand{\z}{\zeta}
\newcommand{\ov}{\overline}
 \renewcommand{\d}{\delta}

\renewcommand\Im{\operatorname{Im}}
\renewcommand\Re{\operatorname{Re}}

\def\qqq{\mathrel{\subset\mkern-15mu\lower.38ex\hbox{${\scriptscriptstyle\rightarrow}$}}}

\let\cal\mathcal

\let\Bbb\mathbb

  \DeclareMathOperator{\spec}{spec}
  
\DeclareMathOperator{\BMO}{BMO}

\DeclareMathOperator{\sing}{sing}
  
\begin{document}

 \title  [A new representation
of  Hankel   operators ] {  A new representation
of  Hankel   operators and its  spectral consequences}



 
\author{ D. R. Yafaev}
\address{ Univ Rennes,   CNRS, IRMAR - UMR 6625, F-35000 Rennes, FRANCE
 and SPGU, Univ. Nab. 7/9, Saint Petersburg, 199034 RUSSIA}
\email{yafaev@univ-rennes1.fr}

\dedicatory{To the memory of Mikhail Zakharovich Solomyak}


\thanks {Supported by  project   Russian Science Foundation   17-11-01126}

\begin{abstract}
We describe a new representation
of  Hankel  operators  $H$ as pseudo-differential operators $A$ in the space of functions defined on the whole axis. The amplitudes of such operators $A$ have a very special structure: they are products of functions of a one variable only. This representation has numerous spectral consequences both for compact Hankel operators and for operators with the   continuous spectrum. 
   \end{abstract}

\maketitle


\section{Introduction }  

{\bf 1.1.}
This  is a short survey based on the talk given by  the  author at the  9th Saint-Petersburg Spectral Theory Conference   held in Euler Institute (Saint-Petersburg, Russia) during 3-6 July 2017. 

Among numerous papers of
M.~Sh.~Birman and M.~Z.~Solomyak on spectral theory of self-adjoint operators, their study (summarized in  \cite{MMS}) of the Weyl asymptotics of eigenvalues of differential operators plays a distinguished role. The methods developed in their papers on this subject were extended by the authors to pseudo-differential and integral operators in \cite{BS4} and \cite{BS3}. We directly use the results of \cite{BS4, BS3} in this article.

\medskip

{\bf 1.2.}
Hankel operators can be defined as integral operators  
\begin{equation}
(H u)(t) = \int_{0}^\infty h(t+s) u(s)ds 
\label{eq:H1}\end{equation}
in the space $L^2 ({\Bbb R}_{+}) $ with kernels $h$ that depend  on the sum of variables only. If necessary, we write $H=H(h)$ for the operator \e{eq:H1}. We refer to the books \cite{NK,Pe} for basic information on Hankel operators. Of course $H$ is symmetric if $  h(t)=\ov{h(t)}$. There are very few cases when Hankel operators can be explicitly diagonalized. The   simplest and most important example $h(t)=t^{-1}$ was considered by T.~Carleman. The corresponding Hankel operator \e{eq:H1} is bounded but not compact; actually, it has the absolutely continuous spectrum $[0,\pi]$ of multiplicity $2$. It follows that a Hankel operator $H$ is compact if, for example, $h\in L^\infty_{\rm loc} ({\Bbb R}_{+})$ and  $h(t)= o(t^{-1})$ as $t\to\infty$ and as $t\to 0$. It turns out that singular values $s_{n} (H)$ of $H$ (and its eigenvalues in the self-adjoint case) have power asymptotics as $n\to \infty$ if the kernel $h(t) $ is close to $t^{-1}$ in the logarithmic scale (both for large and small  $t$), that is, if $h(t)$ behaves as $ \kappa_{\infty }t^{-1}  |\ln t |^{-\alpha}$ with some $\alpha>0$  for $t\to\infty$ and as $ \kappa_0 t^{-1}  |\ln t |^{-\alpha}$ for $t\to 0$.
On the contrary, $H$ is unbounded if    $ h(t)  t \to\infty$ as $t\to\infty$ or as $t\to 0$.

Our first goal is to describe in Section~2 a procedure suggested in \cite{Y, Yafaev3}  reducing an arbitrary  Hankel operator $H$ by an explicit unitary transformation ${\bf M}$ (essentially, by the Mellin transform) to a special integral, or pseudo-differential, operator $A$ in the space $L^2 ({\Bbb R}) $:  
    \begin{equation}
H= {\bf M} ^{-1} A {\bf M}.
 \label{eq:ig6}\end{equation}
  In many cases, the spectral properties of the  operators $A$ are easier to study than those of the original Hankel operators $H$. We emphasize that the identity \e{eq:ig6} does not require that the operators $H$ be symmetric, but in our spectral applications $H$ are self-adjoint (except Section~6).

The operator $A$ can be defined as follows. Put
\[
(X u)(x)= x u(x),\q (D u)(x)=- iu' (x).
\]
Then
  \begin{equation}
A= v(X) s(D) v(X)
 \label{eq:Sam}\end{equation}
where  the standard  function
    \begin{equation}
v(x) = \frac{\sqrt{\pi} }{\sqrt{\cosh(\pi x)} }
   \label{eq:vv}\end{equation}
 is quite   explicit; it is the same for all Hankel operators. The function
  $s(\xi)$ depends of course on $h(t)$, and it can be constructed in the following way. 
  
  Let us formally define the so-called {\it sigma-function}     by the equation 
  \begin{equation}
  h (t)  =    \int_0^\infty     e^{-t \lambda} \sigma(\lambda)  d \lambda  ; 
   \label{eq:conv1}\end{equation}
     in general, $\sigma(\lambda)$ is a distribution.  The function $s(\xi)$ (it is called  the {\it  sign-function} of a Hankel operator $H$ in \cite{Y}) differs from $\sigma(\lambda)$ by a change of variables only:
   \begin{equation}
s(\xi)=\sigma (e^{- \xi}).
 \label{eq:sig2}\end{equation}
 Thus the operator $A$   can be considered  in the space $L^2 ({\Bbb R})$ either as
 a $\Psi$DO (pseudo-differential operator)      with  the amplitude  
  \[
   a(x,y ; \xi)=  v(x) s(\xi) v(y), \q x, y, \xi \in {\Bbb R},
\]
or as an integral operator   with kernel
    \[
(2\pi)^{-1/2}  v(x) (\Phi^*  s )(x-y) v(y) 
 \]
 where $\Phi$,
  $$
 (\Phi  f)(\xi)= (2\pi)^{-1/2}\int_{-\infty}^\infty e^{- ix \xi} f(x) dx,
 $$ 
 is the   Fourier transform. Of course a correct definition of the operator $A$  requires some assumptions on $h(t)$; this will be discussed later.

\medskip

{\bf 1.3.}
We apply this construction to two essentially different classes of  Hankel operators $H$ which lead to two essentially different classes of   $\Psi$DO $A$.

In Section~3, we consider  generalized Carleman operators with kernels  
 \begin{equation}
h(t)= P (\ln t)t^{-1}
\label{eq:LOG}\end{equation} 
where 
 \begin{equation}
P (\xi)=\sum_{m=0}^n p_m \xi^m, \q p_n =1,
\label{eq:LOG4}\end{equation}
 is an arbitrary real polynomial. 
 Obviously, kernels \e{eq:LOG} have two singular points $t=\infty$ and $t=0$.
For $n\geq 1$, such  Hankel operators   are unbounded but are well defined as self-adjoint operators.

For kernels \e{eq:LOG}, the sign-function $s(\xi) $ is  a  real polynomial  
 \begin{equation}
 Q(\xi)=  \sum_{ m=0}^n q_m \xi^m 
\label{eq:LO4}\end{equation}
  determined by $P (\xi )$. In this case
  \begin{equation}
A= v (X) Q (D ) v (X)  
\label{eq:LOGz}\end{equation}
 is a differential operator.
  The polynomials $P (\xi)$ and $Q (\xi )$ have the same degree,   and their coefficients   are linked by an explicit formula (see formula \e{eq:LOG5a} below); in particular, $ q_{n}   = (-1)^n$. If $n=0$, then $Q(\xi )=P( \xi)=1$ so that $A$ is the multiplication operator by $v(x)^2$. This yields the familiar diagonalization of the Carleman operator.

    Observe that the highest order term of the operator $A$ equals $v^2 (x)  D ^n $
 where $  v^2 (x) $ tends to zero (exponentially) as $| x |\to\infty$.
                    Apparently such differential operators  were never studied before, and we are led to fill in this gap.
Studying   differential operators \e{eq:LOGz} in Section~3,  we do not make specific assumption \e{eq:vv} and consider sufficiently arbitrary real  functions $v(x)$  tending to zero as $| x | \to\infty$.
     The essential spectrum of differential operators \e{eq:LOGz} 
 was  localized in \cite{Yf1}    where it was shown that $\spec_{\rm ess} (A)=\spec  (A)={\Bbb R}$ if $n$ is odd, and $\spec_{\rm ess} (A)= [0, \infty)$ if $n$ is even. 
   The last result should be compared with the fact that $\spec_{\rm ess} (A)= [\min Q( \xi ), \infty)$  if $v(x)=1$. Thus, even in this relatively simple question,  the degeneracy of $v(x )$ at infinity significantly changes spectral properties of  differential operators $A$.
       The detailed spectral structure, in particular, the absolutely continuous spectrum,  of differential operators \e{eq:LOGz} and hence of the Hankel operators with kernels \e{eq:LOG} was described in \cite{Ya}.

\medskip

{\bf 1.4.} 
In Section~4, we are interested in compact self-adjoint Hankel operators \e{eq:H1} with power-like asymptotics of eigenvalues $\lambda_n^\pm(H)$ as $n\to\infty$. Let us denote by $\{\lambda_n^+(H)\}_{n=1}^\infty$  the non-increasing sequence
of positive eigenvalues of a compact 
self-adjoint operator $H $ (with multiplicities taken into account), and set $\lambda_n^-(H)=\lambda_n^+(-H)$.
Sharp estimates of $\lambda_n^\pm(H)$ (and, more generally, of singular values $s_n (H)=\lambda_n^+(\sqrt{H^* H})$ in the non-self-adjoint case) are very well known. Thus V.~V.~Peller found (see Chapter~6 of his book \cite{Pe}) necessary and sufficient conditions for the validity of the estimates $s_n (H)=O (n^{-\alpha})$.
At the same time, there are practically no results on the asymptotic behaviour
of   eigenvalues of Hankel operators. 
This state of affairs is in a sharp contrast with the case of differential operators, 
where the Weyl type asymptotics of eigenvalues is established in a very large variety of situations. 
Our goal here is to fill in this gap by 
describing  classes of Hankel operators where the leading term of eigenvalue asymptotics  
can be found explicitly.

In general, the study of eigenvalue asymptotics for any class of operators
involves two steps: 
construction of an appropriate model problem where the eigenvalue asymptotics
can be determined sufficiently explicitly, and using eigenvalue estimates  (or variational methods)
to extend the asymptotics   to a wider class of operators. Apparently, for a given Hankel operator $H$, there is no natural model operator in the class of Hankel operators. 
So, the crucial step of our approach is a construction of the model operator $\Psi$DO  $A_{*}$. The spectral asymptotics of $A_{*}$ and hence of the corresponding  Hankel operators $H_{*}$ is given by the Weyl law. This result is then extended to the initial operator $H$.

 \medskip

{\bf 1.5.} 
Hankel operators  can also be naturally realized in the space $\ell^2({\Bbb Z}_{+})$ of sequences. Namely,
for  sequences $g=\{g (j) \}_{n\in {\Bbb Z}_{+}}$,
Hankel operators $G=G (g)$ are defined  by   infinite matrices:
\begin{equation} 
(G u) (j)=\sum_{k=0}^\infty g (j+k) u (k).
 \label{eq:Hd}\end{equation}
 As is well known, the operators $H$ and $G$ give two different representations of the same object.
 Indeed, let us introduce a  unitary transformation 
  $U: \ell^2({\Bbb Z}_{+}) \to L^2({\Bbb R}_{+})$ by the formula 
     \[
(U u)(t)=  \sum_{j=0}^\infty {\cal L}_{j } (t)   u (j) e^{-t/2},\q u=\{ u (j) \}_{j \in {\Bbb Z}_{+}}\in \ell^2({\Bbb Z}_{+})
\]
 where ${\cal L}_{j}={\cal L}_{j}^0$ are the Laguerre polynomials.
   The Hankel operators $H$ and $G$ are linked by this transformation, that is
   \begin{equation}
U G U^*=H, 
 \label{eq:LL1}\end{equation}
 provided
  \[
h(t)= \sum_{j=0}^\infty {\cal L}_{j}^1 (t)   g (j) e^{-t/2}
\]
 where $ {\cal L}_{j}^1$ are the generalized Laguerre polynomials.
 
 Although the operators $H$ and $G$ are unitary equivalent, it is convenient to study eigenvalues asymptotics of $G$ independently of the results on $H$. This is done in Section~5.
 
 In the discrete case the role of the Carleman operator is played by the Hankel operator $G$ (the Hilbert matrix) with the matrix elements $g(j)=(j+1)^{-1}$. This operator is bounded but not compact; actually, it has the simple absolutely continuous spectrum $[0,\pi]$. Note that the asymptotic behavior of a kernel $h(t)$ like $ t^{-1}  |\ln t |^{-\alpha}$ as $t\to\infty$ (resp. as $t\to 0$) is essentially equivalent to the asymptotic behavior of  the  matrix elements of  the corresponding Hankel operator $G$ like $ (-1)^j j^{-1}  (\ln j )^{-\alpha}$ (resp. like $   j^{-1}  (\ln j)^{-\alpha}$) as $j\to \infty$. It is often useful to keep in mind that Hankel operators $G$ and $\wt{G}$ are unitarily equivalent provided their matrix elements are linked by the relation $\wt{g} (j)=(-1)^j g(j)$.
 
It is natural to expect that a faster rate of convergence to zero
as $j\to\infty$
 of the sequence $g(j)$  results in a faster convergence   to zero
as $n\to\infty$ of the eigenvalues $\lambda^\pm_n(G)$. 
Indeed, there is a  deep result of H.~Widom who showed  in \cite{Widom} that for $\gamma>1$ 
the Hankel operator corresponding to the sequence $g(j)=(j+1)^{-\gamma}$  is non-negative   and its eigenvalues converge to zero as
$$
\lambda_n^+(G)=\exp(-\pi \sqrt{2\gamma n}+o(\sqrt{n})), \quad n\to\infty.
$$
By the way, for the proof of this result, H.~Widom also used a reduction of Hankel operators he considered to $\Psi$DO $
s(D)^{1/2}v(X)^2 s(D)^{1/2}$. Such a reduction  is possible if $s(\xi)\geq 0$.

    Note that both Hankel operators $H$ and $G$ can be realized in Hardy spaces of analytic functions where these operators are determined by their symbols. We do not discuss the representations of Hankel operators in Hardy spaces but emphasize that their symbols briefly mentioned in Section~6 and sigma-functions are completely different objects --
    see Section~3 in \cite{Yafaev3}.

 \medskip

{\bf 1.6.}
    Finally, in Section~6 we discuss more general results on singular values and eigenvalues of Hankel operators with kernels $h(t)$ oscillating as $t\to \infty$ and with matrix elements $g(j)$ oscillating as $j\to \infty$.
    
    The results on power-like asymptotics of singular values  have direct applications to rational approximations of functions with logarithmic singularities. A result of such type is stated in Section~6.

\section{Main identity }  

A detailed presentation of
the results of this section can be found in the papers \cite{Y, Yafaev3}.

\medskip

{\bf 2.1.}
For a given Hankel operator  $H$,
 let the sigma-function $\sigma(\lambda)$ be  {\it formally}  defined by equation \e{eq:conv1}, and let $\Sigma$ be the operator of multiplication by $\sigma$, that is,
   \begin{equation}
  (\Sigma  f) (\lambda)= \sigma(\lambda) f(\lambda),\q \lambda>0.
\label{eq:Ssi}\end{equation} 
We will show that
\begin{equation} 
H= {\sf L}^* \Sigma {\sf L} 
 \label{eq:MAIDp}\end{equation} 
 where  ${\sf L}$ is the Laplace transform:
  \begin{equation} 
  ({\sf L} u) (\lambda)= \int_{0}^\infty e^{-t \lambda} u(t) dt  .
\label{eq:LAPj}\end{equation}  
A {\it formal} proof   of  the identity \e{eq:MAIDp}
 is quite simple. Indeed, the integral kernel of the operator in the right-hand side of  \e{eq:MAIDp} equals
 \[
 \int_0^\infty e^{-\lambda t} \sigma (\lambda) e^{-\lambda s} d \lambda=h(t+s)
 \]
 if $ \sigma (\lambda)$ and $h(t)$ are linked by formula  \e{eq:conv1}. Thus it equals the integral kernel of the  operator  defined by \e{eq:H1}.

 The precise sense of formula \e{eq:MAIDp} needs of course  to be clarified. Observe that,
 by its definition \e{eq:conv1},  $\sigma(\lambda)$  can be a regular function only for   kernels $h(t)$      satisfying  some specific analytic  assumptions. Without such very restrictive assumptions,  $\sigma$ is necessarily a distribution. Even for very good kernels $h(t)$ (and especially for them), $\sigma(\lambda)$ may be a highly singular distribution.
For example,    for $h(t)=t^k e^{-\alpha t}$ where  $\Re \alpha >0$ ($\alpha$ may be complex) and $k=0,1, \ldots$, the sigma-function $\sigma(\lambda)= \d^{(k)}(\lambda-\alpha )$ is a  derivative   of the delta-function. On the contrary, singular kernels $h(t)$ may yield    sigma-functions $\sigma(\lambda)$   smooth on ${\Bbb R}_{+}$.  For example, if $h(t)=t^{-q}$ where $q>0$ may be arbitrary large, then $\sigma(\lambda)= \Gamma(q)^{-1}\lambda^{1-q}$; here and below $\Gamma(\cdot)$ is the gamma function. 

Thus, we replace \e{eq:MAIDp} by the identity    
  \[
 (H f_{1}, f_{2}) =  ( \Sigma {\sf L} f_{1}, {\sf L} f_{2})
\] 
for arbitrary test functions $f_{1}, f_{2} \in  C_{0}^\infty ({\Bbb R}_{+})$.  With respect to $h$ we require only that  $h \in  C_{0}^\infty ({\Bbb R}_{+})'$.  Let us introduce the Laplace convolution
 \[
( \bar{f}_{1}\star f_{2} )(t)=\int_{0}^\infty \ov{f_{1}(s)} f_{2} (t-s)ds.
 \]
 Then formally
 \[
 (H f_{1}, f_{2})=\la h,  \bar{f}_{1}\star f_{2}\ra.
  \]

By one of the versions of the Paley-Wiener theorem, the Laplace transform $\sf L$ is an isomorphism of $C_{0}^\infty ({\Bbb R}_{+})$ onto the space ${\cal Y}$   of analytic functions $g(\lambda)$ of $\lambda\in{\Bbb C}$ exponentially decaying as $\Re \lambda \to +\infty$, exponentially bounded as $\Re \lambda \to -\infty$ and decaying faster than any power of
  $|\lambda|^{-1}$ as $|\Im \lambda| \to \infty$ (see \cite{Yafaev3} for details).   By duality, ${\sf L}^*: {\cal Y}' \to C_{0}^\infty    ({\Bbb R}_{+})'$ is also an isomorphism, and  hence according to definition \e{eq:conv1}, 
   \begin{equation} 
  \sigma= ({\sf L}^*)^{-1}h \in {\cal Y}'
  \label{eq:ss}\end{equation}
if  $h \in  C_{0}^\infty ({\Bbb R}_{+})'$. This yields a one-to-one correspondence between kernels $h\in C_{0}^\infty ({\Bbb R}_{+})'$ of Hankel operators and their sigma-functions $\sigma  \in {\cal Y}'$ and makes the theory   self-consistent. 
Note that instead of operators, we consequently work with quadratic forms  which is both more general and more convenient. For $g\in {\cal Y}$, we set $g^*(\lambda)=
  \ov{g(\bar{\lambda})} $.

 Now we are in a position to   state our main  identity \e{eq:MAIDp} precisely. 
  
 \begin{theorem}\label{1}
Let   $h \in  C_{0}^\infty ({\Bbb R}_{+})'$, and let $\sigma \in {\cal Y}'$ be defined by formula
  \e{eq:ss}.
  Then       the identity
 \[
 { \pmb\la} h ,  \bar{f}_{1}\star f_{2} {\pmb \ra} =   { \la} \sigma , (  {\sf L} f_{1} )^*{\sf L} f_{2} { \ra}
\]
  holds for arbitrary
   $f_{1}, f_{2}\in C_{0}^\infty    ({\Bbb R}_{+})$.
 \end{theorem}
 
 The identity \e{eq:MAIDp} does not of course give a diagonalization of  Hankel operators because
the operator ${\sf L} $ is not unitary. However it is  continuously invertible as a mapping ${\sf L} : C_{0}^\infty    ({\Bbb R}_{+}) \to {\cal Y}$ so that equality  \e{eq:MAIDp} plays the same role as Sylvester's inertia theorem  which states that two Hermitian matrices $H$ and $\Sigma$ related by equation \e{eq:MAIDp} have the same  total numbers of   positive   and negative   eigenvalues. In particular,      $\pm H\geq 0$ if and only if $\pm \Sigma\geq 0$. According to Theorem~\ref{1} the same assertion is true for Hankel operators $H$ and operators of multiplication $\Sigma$. Now the operators $H$ and $\Sigma$ are of a completely different nature and $\Sigma$ (but not $H$) admits an explicit spectral analysis. 
As an example of this approach, we show in Section~4 of  \cite{Yf1a} that if $\sigma (\lambda) >0$ (or $\sigma (\lambda) <0$) on a set of positive Lebesgue measure, then the Hankel operator $H$ has infinite positive (or negative) spectrum. On the other hand, singularities of $\sigma(\lambda)$ at some isolated points produce   finite numbers (depending on the order of the singularity) of positive or negative eigenvalues (see Section~4 of  \cite{Y}). In particular, this approach enables us \cite{Yf} to give an explicit formula for total numbers of positive and negative eigenvalues of finite-rank Hankel operators.

\medskip
 
 {\bf 2.2.}
 To perform the spectral analysis of Hankel operators, we will transform the identity \e{eq:MAIDp} using the factorization of the   operator $\sf L$. Let us introduce the Mellin transform $M :L^2 ({\Bbb R}_{+}) \to L^2 ({\Bbb R} ) $,
 \[
( M  u)(x)=  (2\pi)^{-1/2} \int_{0}^\infty u(t) t^{-1/2 - i x} dt , 
 \]
 the reflection operator ${\cal J}$, $({\cal J}u)(x)= u(- x)$,  
 and set
   \[
  ({\pmb \Gamma}    u)(x)=\Gamma(1/2 +i x) u(x). 
 \]
 We use the following elementary fact.

   \begin{lemma}\label{LAPL} 
  For the Laplace transform defined by \e{eq:LAPj}
    the identity
    \begin{equation}
({\sf L}f )(\lambda)=  (M ^{-1} {\cal J} {\pmb \Gamma}   M f )(\lambda), \q \lambda>0,
 \label{eq:LAPL1}\end{equation}
 holds for all $f\in L^2 ({\Bbb R}_{+})$.
  \end{lemma}   
  
By the way,    factorization  \e{eq:LAPL1} enables one to invert the Laplace transform:
   \[
{\sf L}^{-1}=  M ^{-1}   {\pmb \Gamma} ^{-1} {\cal J}  M   .
 \]
 We use \e{eq:LAPL1} to establish the unitary equivalence of the operators $H$ and $A$. 
 Observe that
 \[
 |\Gamma (1/2+ i x)|= \frac{\sqrt{\pi} }{\sqrt{\cosh(\pi x)} }= v( x)
\]
and set
 \[
 ({\bf M} f) (x)= \frac{\Gamma (1/2- i x)}{|\Gamma (1/2- i x)|} (M  f) ( -x).
\]
Then   \e{eq:LAPL1} can be rewritten as
   \begin{equation}
{\sf L}= M ^{-1} v(X) {\bf M}.
  \label{eq:SX}\end{equation}
 
 Note also that $M=\Phi^{-1} W$ where the unitary operator $W:L^2 ({\Bbb R}_{+})\to L^2 ({\Bbb R} )$ is defined by $(W f)(\xi)= e^{-\xi/2} f(e^{-\xi})$, and hence it follows from  \e{eq:sig2},  \e{eq:Ssi}  that
   \[
s(D)=M^{-1}\Sigma M.
\]
  In view of \e{eq:SX} this yields the identity 
   \begin{equation}
H= {\bf M} ^{-1} v(X) s(D) v(X) {\bf M}.
 \label{eq:sig6}\end{equation}

 Let us state a simple sufficient condition for the validity of this representation which has been derived above rather formally.

   \begin{theorem}\label{SIG} 
   Suppose that 
   \[
\sigma \in L^\infty ({\Bbb R}_{+})
\]
 so that   the   operator $s(D) $ is bounded.
  Then the identity
 \e{eq:sig6} holds.
   \end{theorem}

\medskip
 
 {\bf 2.3.}
 The above results on integral operators \e{eq:H1} can be extended to Hankel operators $G$ defined by formula \e{eq:Hd} in the space $\ell^2({\Bbb Z}_{+})$.
  In the discrete case, the role   of the sigma-function $\sigma (\lambda)$ of $\lambda\in{\Bbb R}_{+}$ is played by  the function $\eta(\mu)$ defined on the interval $(-1,1)$ and linked to $\sigma (\lambda)$  by the relation
   \begin{equation}
\sigma (\lambda)=\eta\left(\frac{2\lambda-1}{2 \lambda+1}\right) .
 \label{eq:DC1}\end{equation}
  Let
  \begin{equation}
g(j)=\int_{-1}^1 \eta(\mu)\mu^j d\mu, 
\quad
j=0,1,2,\dots,
 \label{eq:DC}\end{equation}
be the sequence of moments of $\eta$, and let   $G$ be the Hankel operator    with the matrix elements $g(j)$. We emphasize that equations \e{eq:DC} play the role of \e{eq:conv1}.
 It can be easily shown that relation \e{eq:LL1} is satisfied if  the kernel $h(t)$ of $H$ is given by  \e{eq:conv1}.
 Therefore an analogue of Theorem~\ref{SIG} is stated as follows.

\begin{theorem}\label{thm.e2}\cite[Theorem 7.7]{Yafaev3}
Let $\eta\in L^\infty(-1,1)$.
Then the Hankel operator $G$ in $\ell^2(\Bbb Z_+)$ with the matrix elements \e{eq:DC}   is unitarily equivalent to the $\Psi$DO 
\eqref{eq:Sam} in $L^2(\Bbb R) $ with the sign-function
$s(\xi) $ defined by \e{eq:sig2}, \e{eq:DC1}.
\end{theorem}

\section{Generalized Carleman operators }

A detailed presentation of
the results of this section can be found in the papers \cite{Yf1, Ya}.

\medskip

{\bf 3.1.} 
Our goal      here is to study spectral properties of generalized Carleman operators with kernels  
\e{eq:LOG} 
where $P (\xi)$
 is an arbitrary real polynomial \e{eq:LOG4}.   In this case, we have
\begin{equation} 
h(t)=\int_{0}^\infty e^{-\lambda t} \sum_{m=0}^n q_{m} \ln^m \lambda d\lambda
\label{eq:MM1}\end{equation} 
where the coefficients
 \begin{equation}
q_m=  (-1)^m  \sum_{j =m}^n   \tbinom{j}{m} \gamma^{(j-m)} (0) p_j ,\q m=0,\ldots, n, \q \gamma ( z)=\Gamma (1-z)^{-1}.
\label{eq:LOG5a}\end{equation} 
For example,
$(-1)^n q_{n}=   p_n $ and
 $ (-1)^n q_{n-1}= -   p_{n-1}  +   \Gamma' (1) \,n\,  p_n  $
    for all $n$ (recall that  $  -\Gamma' (1)$ is the Euler constant).  Of course formulas \e{eq:LOG5a} enable one to recover the coefficients $p_{n}, p_{n-1},\ldots, p_{0}$ given the coefficients $q_{n}, q_{n-1},\ldots, q_0$. It follows from \e{eq:MM1} that
$s(D)=: Q(D)$ is a differential operator given by formula  \e{eq:LO4}.

For Hankel operators \e{eq:H1} with kernels  \e{eq:LOG4}, the identity \e{eq:sig6} yields the following result.

 \begin{theorem}\label{1z} \cite[Theorem~3.2]{Yf1}
 Let  $Q(\xi)$ be polynomial \e{eq:LO4} with the coefficients $q_{m}$   defined by formulas \e{eq:LOG5a}, and let $A  $
 be the differential operator \e{eq:LOGz}. Then for all functions $u_{j} $, $j=1,2$, such that their Mellin transforms
 $M u_{j} $ belong to $ C_{0}^\infty ({\Bbb R})$,
 the identity 
 \[
(H u_{1}, u_{2})=  (A {\bf M} u_{1}, {\bf M} u_{2})
 \]
holds.
 \end{theorem}

A large part of our results on generalized Carleman operators can be summarized by the following assertion.     Below we denote by $\la x\ra$ the operator of multiplication by the function $(1+x^2)^{1/2}$.

  \begin{theorem}\label{SpThHa}
Let $H $ be the self-adjoint Hankel   operator defined by formula \e{eq:H1} where $h(t)$ is function  \e{eq:LOG}  and $P(\xi)$ is a  real polynomial  \e{eq:LOG4} of degree $n\geq 1$. 
Then
\begin{enumerate}[\rm(i)]
\item
The spectrum of the operator $H$ is absolutely continuous except eigenvalues that  may accumulate to zero and  infinity only.   
   \item
The 
 absolutely continuous  spectrum of the operator $H$ covers $\Bbb R$ and is simple for $n$ odd. It coincides with $[0,\infty)$ and has multiplicity $2$  for $n$ even. 
 \item
 If $n$ is odd, then the multiplicities of eigenvalues  of the operator $H$  are bounded $($from above$)$ by $(n-1)/2$. If $n$ is even, then the multiplicities of positive eigenvalues are bounded by $n/2-1$, and the multiplicities of negative eigenvalues are bounded by $n/2$. 
  \item
 For any $\d>1/2$, the operator-valued  function 
      \begin{equation}
  \la \ln t\ra^{-\d} ( H-z)^{-1} \la \ln t\ra^{-\d}, \q \Im z \neq 0,
 \label{eq:LAPH}\end{equation}
 is H\"older continuous with exponent $\alpha<\d-1/2$ $($and $\alpha<1)$ up to the real axis,  except the eigenvalues of the operator $H$ and the point zero.
\end{enumerate}
    \end{theorem}

    Clearly, this assertion is similar in spirit to the corresponding results for differential operators of Schr\"odinger type. 
    Statement (iv) is known as  the limiting absorption principle.

 \medskip
    
     {\bf 3.2.}
        In view of Theorem~\ref{1z}, the proof of Theorem~\ref{SpThHa} reduces to a proof of the corresponding results for the differential operator $A$ defined by \e{eq:LOGz}. However the standard results on differential operators are not applicable in this case because of
                a strong degeneracy of $v(x)$ at infinity. Fortunately,
  operators  \e{eq:LOGz}
    can be reduced by an explicit unitary transformation $  {\bf L}$ (the generalized  Liouville transformation)  to standard differential operators. Set
   \[ 
 (  {\bf L} u) (x)=  y'(x)^{1/2} u( y(x))
\]
 where the variables $x$ and $y$ are linked by  the relation
\[
y= y(x)= \int_0^x v( s )^{-2/n} d s
\]
 so that $y'(x)=v(x)^{-2/n}$. Then    
$B= {\bf L}^{-1} A  {\bf L}$
 is also a differential operator in the space $L^2 ({\Bbb R})$, and it is given by the formula    
   \begin{equation}
B=    D^n+ \sum_{m=0}^{n-1} b_m (y) D^m, \q  \q D   = D_{y}   =- id/d y.
 \label{eq:VB}\end{equation}
 Our crucial observation is that   the  coefficients $b_m (y)$, $m=0,1,\ldots, n-1$, of the operator  $B$  decay at infinity. Moreover, all coefficients $b_{0}(y), \ldots, b_{n-2} (y)$  of the operator $B$ are short-range, that is, they decay faster than $|y|^{-1}$ as $|y|\to\infty$. The coefficient $b_{n-1} (y)$   can be removed by a gauge transformation $\cal J$ defined by 
 \[
 ({\cal J}u)(y)=e^{i \phi(y)} u(y)  \q {\rm  where}\q  \phi(y)=-\frac{1}{n}\int_{0}^y b_{n-1} (s) ds.
 \]
 This means that  the operator $\wt{B}= {\cal J}^* B{\cal J}$ has again the form \e{eq:VB} with $\wt{b}_{n-1} (y)=0$. The
 coefficients $\wt{b}_{0}(y), \ldots, \wt{b}_{n-2} (y)$ remain short-range.
 
Using fairly standard methods of scattering theory    we obtain assertions (i), (ii) and (iii) of Theorem~\ref{SpThHa} for the operator $B$. Since 
    \begin{equation}
    A= {\bf L} B  {\bf L}^{-1 } \q {\rm and}\q
H={\bf M}^{-1 }  A {\bf M} ,
   \label{eq:equiv}\end{equation}
these results   remain true  for the operators $A$ and $H$.   
Recall that for differential operators $B$,  one defines their (generalized) eigenfunctions as special solutions  $\psi (y, k)$, $k\in{\Bbb R}$, of the equation $B \psi=k^n\psi$ satisfying some asymptotic conditions as $y\to\infty$ and $y\to -\infty$ and then establishes the expansion theorem over these eigenfunctions. Relation \e{eq:equiv} allows us to  carry over these results to the operators $A$ and $H$. According to \e{eq:equiv} one can define the eigenfunctions of the operator $H$   by the equality
 \begin{multline*}
\theta  (t,k)=  ({\bf M}^{-1 }  {\bf L} \psi ( k)  )(t)
\\
= (2\pi)^{-1/2} t^{-1/2 } \int_{-\infty}^\infty e^{  -i x(y)\ln t}
e^{i \eta(x (y))} x'(y)^{1/2}   \psi  (y,k) d y
\end{multline*}
 where $\eta(x )=\arg\Gamma ( 1/2+ix)$. This integral converges although not absolutely. Applying the stationary phase method one can deduce  asymptotics of the eigenfunctions $\theta  (t,k)$ as $t\to 0$ and as $t\to\infty$  from this representation. In particular, we obtain
 the uniform in $k$ (on compact subsets of ${\Bbb R}\setminus\{0\}$, away from the eigenvalues  of $H$) estimate
     \[
    | \theta(t ,k)|\leq C t^{-1/2}
    \]
     and a similar estimate on  differences $\theta(t ,k')-\theta(t ,k)$. Using these estimates one can prove  assertion  (iv)
     of Theorem~\ref{SpThHa}. This result looks similar to  the limiting absorption principle for differential operators.
   The difference, however, is that the weight  is $  \la \ln t\ra^{-\d} $ in  \e{eq:LAPH} while it is $  \la   x\ra^{-\d} $ (also with $\d>1/2$) for the resolvents of differential operators.
    Thus the power scale for differential operators corresponds to the logarithmic scale for Hankel operators.

 \medskip
    
     {\bf 3.3.}
     Actually, the specific expression \e{eq:vv} for the function $v(x)$ in  definition \e{eq:LOGz} of the operator $A$ is inessential.  It is noteworthy that if $v(x)$    tends to zero exponentially as $|x|\to\infty$, then the coefficients   $b_{0}(y), \ldots, b_{n-2} (y)$ of the operator $B$ decay faster than $|y|^{-1}$ as $|y|\to\infty$. On the contrary, for slower decay of $v(x)$, these   coefficients     decay slower than (or as)   $|y|^{-1}$.     Thus, somewhat counter-intuitively,  a stronger     degeneracy of the operator \e{eq:LOGz}  yields  better properties of the operator $B= {\bf L}^{-1} A  {\bf L}$.
     
    We also note that our
   approach applies to  sufficiently arbitrary differential operators of order $n$  with a degeneracy of the coefficient in front of $D^n$.

\section{Compact operators. Asymptotics of singular values and eigenvalues}

 {\bf 4.1.  } 
In this subsection we collect necessary auxiliary results.
First we recall Weyl asymptotics of $\Psi$DO. For
  $x\in {\Bbb R}$, we  use the standard notation $x_{\pm}=\max\{0,\pm x\}$.
 
\begin{theorem}\label{thm.b2}
Let $s\in C^\infty ({\Bbb R})$ be a  real-valued 
function such that
\begin{equation}
s(\xi)=
\begin{cases}
s_{\infty} \xi^{-\alpha}(1+o(1)), &
\xi\to\infty,
\\
s_{-\infty} |\xi|^{-\alpha}(1+o(1)), &
\xi\to-\infty,
\end{cases}
\label{b2}
\end{equation}
for some   $\alpha>0$ and some constants  $s_{\infty}$ and $s_{-\infty}$.
Assume that $v(x)=\overline{v(x)}$ and
\begin{equation}
| v(x)| \leq C\la x\ra ^{-\rho}, 
\quad 
x\in\Bbb R,
\label{b1a}
\end{equation}
for some $\rho>\alpha/2$. Put
\[
{\bf a}^{\pm}
=
(2\pi)^{-\alpha}
\bigl(
(s_{-\infty})_\pm^{1/\alpha}+ (s_{\infty})_\pm^{1/\alpha}
\bigr)^\alpha
\biggl(
\int_{-\infty}^\infty  | v(x)| ^{2/\alpha}dx\biggr)^\alpha.
\]
Then for the $\Psi$DO \eqref{eq:Sam}  in $L^2(\Bbb R)$ one has 
\[
\lambda_n^\pm(A)
=
{\bf a}^{\pm} n^{-\alpha}+o(n^{-\alpha}), 
\quad 
n\to\infty.
\]
\end{theorem}

For compactly supported $v$, 
Theorem~\ref{thm.b2}    
was proven by M.~Sh.~Birman and M.~Z.~Solomyak in \cite{BS4}; actually,  the multi-dimensional case was considered
 in \cite{BS4}. Their result can be easily extended   to arbitrary functions $v$ satisfying \eqref{b1a} (see, e.g., Appendix in \cite{II}).
For general $\Psi$DO (acting in a bounded domain) with amplitudes asymptotically homogeneous at 
infinity, Weyl type formula for the asymptotics of the  spectrum was obtained in \cite{BS3}. 

We need also some estimates on singular values $s_{n} (H) $ of Hankel operators \e{eq:H1}. They  are stated in the next  assertion established in \cite{I}. For $\alpha\geq 1/2$, the  proof of these estimates relies heavily on deep results by V.~V.~Peller (see Chapter~6 of his book \cite{Pe}).
For an arbitrary $\alpha>0$, we set 
\begin{equation} 
 N(\alpha)=[\alpha]+1\; {\rm if} \; \alpha\geq 1/2 \q {\rm and } \q N(\alpha) =0\; {\rm if} \;\alpha<1/2.
\label{eq:Mm}\end{equation}

\begin{theorem}\label{thm.a2} 
Let $\alpha>0$ and  let  $h \in L^\infty_{\rm loc}({\Bbb R}_+)$ be a complex valued function; if $\alpha\geq 1/2$, suppose also 
that $h\in C^{N(\alpha)}({\Bbb R}_+)$.
Assume that $h$ satisfies the conditions
$$
h^{(m)}(t)=o(t^{-1-m}\la\log t\ra^{-\alpha})
\quad 
\text{ as $t\to 0$ and as $t\to\infty$}
$$
for all $m=0,1,\ldots, N(\alpha)$.
Then 
$s_n(H)=o(n^{-\alpha})$ as $ n\to\infty$.
\end{theorem}

Theorem~\ref{thm.a2} remains true if   $o (\cdot)$ (on both occasions above) is replaced by  $O (\cdot)$.

 We need also
the following standard result (see, e.g.,  \cite[Section 11.6]{BSbook})
in spectral perturbation theory, which asserts the stability of eigenvalue asymptotics.

\begin{lemma}\label{lma.b1} 
Let $K_{0}$ and $K$ be compact self-adjoint operators and let $\alpha>0$.
Suppose that, for both signs $``\pm"$,   
$$
\lambda^\pm_n(K_{0})
=
{\bf a}^\pm n^{-\alpha}+o(n^{-\alpha}),
\quad \text{ and }\quad
s_n(K)
=
o(n^{-\alpha}),
\quad
\quad n\to\infty.
$$
Then 
$$
\lambda^\pm_n(K_{0}+K)
=
{\bf a}^\pm n^{-\alpha}+o(n^{-\alpha}),
\quad n\to\infty.
$$
\end{lemma}

    \medskip

  {\bf 4.2.}
  Our main result on asymptotics of eigenvalues of Hankel operators \e{eq:H1} can be stated as follows.
Let us set
  \begin{equation}
\tau(\alpha) =2^{-\alpha}
\pi^{1-2\alpha}
B(\tfrac1{2\alpha},\tfrac12)^{\alpha}
\label{a4}
\end{equation}
where 
\[
B(x,y) =
\frac{\Gamma(x)\Gamma(y)}{\Gamma(x+y)}
\]
 is the 
Beta function.

\begin{theorem}\label{thm.d1}
Let $\alpha>0$ and let the integer $N(\alpha)$ be given by \e{eq:Mm}.
Let $h$ be a real valued function in $L^\infty_{\rm loc}({\Bbb R}_+)$; if $\alpha\geq1/2$, assume also 
that $h\in C^{N(\alpha)}({\Bbb R} _+)$.
Suppose that 
\begin{align*}
\biggl( \frac{d}{dt}\biggr)^m
\bigl(h(t)-\kappa_0 t^{-1}(\log(1/t))^{-\alpha}\bigr)
&=
o(t^{-1-m}\la\log t\ra^{-\alpha}), \quad t\to0, 
\\
\biggl( \frac{d}{dt}\biggr)^m
\bigl(h(t)-\kappa_\infty t^{-1}(\log t)^{-\alpha} \bigr)
&=
o(t^{-1-m}\la\log t\ra^{-\alpha}), \quad t\to\infty,
\end{align*}
for some $\kappa_0, \kappa_\infty\in\Bbb R$ and all $m=0,\dots, N(\alpha)$. 
Then the eigenvalues of the corresponding Hankel operator  $H$ 
have the asymptotic behaviour
\begin{equation}
\lambda_n^\pm(H)
=
a^\pm n^{-\alpha}
+
o(n^{-\alpha}) ,\q n\to\infty,
\label{d3a}
\end{equation}
where
\begin{equation}
a^\pm
=
\tau (\alpha)\bigl((\kappa_0)_\pm^{1/\alpha}+(\kappa_\infty)_\pm^{1/\alpha}\bigr)^\alpha.
\label{eq:Al}
\end{equation} 
\end{theorem}


Our proof of this result relies on the following three ingredients:

\begin{enumerate}[(i)]
\item
Theorem~\ref{SIG} allows us to replace
 Hankel operators  \e{eq:H1}  by the 
  $\Psi$DO  $A$ defined by \e{eq:Sam}.
\item
 Weyl type  spectral asymptotics for    $\Psi$DO of this type stated in Theorem~\ref{thm.b2}.
\item
Estimates on singular values of Hankel operators of Theorem~\ref{thm.a2}.
\end{enumerate}

  \medskip
    
 {\bf 4.3.  } 
 Let us sketch the proof of Theorem~\ref{thm.d1}. The first and the most important step  
is to construct a \emph{model operator}. To  that end, we introduce an auxiliary explicit function
  by the formula
\begin{equation}
\sigma_*(\lambda)
=
\kappa_\infty\ |\log \lambda|^{-\alpha}\chi_0(\lambda)
+
\kappa_0 |\log \lambda|^{-\alpha}\chi_\infty(\lambda), 
\quad 
\lambda>0,
\label{d6}
\end{equation}
where the    
cut-off functions  $\chi_0,\chi_\infty\in C^\infty(\Bbb R_+)$ satisfy
\begin{equation}
\chi_0(t)=
\begin{cases}
1& \text{for $t\leq1/4$,}
\\
0& \text{for  $t\geq1/2$,}
\end{cases}
\qquad
\chi_\infty(t)=
\begin{cases}
0& \text{for $t\leq2$,}
\\
1& \text{for $t\geq4$.}
\end{cases}
\label{a7a}
\end{equation}

It turns out that the functions $h(t)$ and $h_{*} (t)= ({\sf L}  \sigma_*)(t)$ have the same   asymptotics  as $t\to \infty$;  a similar relation holds as $t\to 0$. For the proof,
we need an  elementary technical result about the Laplace transform of 
functions with logarithmic singularities at $\lambda=0$ and $\lambda= \infty$.

\begin{lemma}\label{L}
Let $\alpha>0$,   $m\in\Bbb Z_+$,
$$
I_m^{(\infty)}(t)= \int_{0}^c (-\log \lambda)^{-\alpha}\lambda^{m} e^{-\lambda t}d\lambda, \q c\in (0,1),
$$
and 
\[
I_m^{(0)}(t)= 
\int_{c}^\infty (\log \lambda)^{-\alpha}\lambda^{m} e^{-\lambda t}d\lambda , \q c>1.
\]
Then 
\begin{equation}
I_m^{(\infty)}(t)=m!\, t^{-1-m} |\log t|^{-\alpha} (1+O(|\log t|^{-1})) 
\label{eq:L2}
\end{equation}
as $t\to\infty$,  and $I_m^{(0)} (t)$ has the same asymptotic behaviour 
\eqref{eq:L2} as $t\to 0$.
\end{lemma}

This result   is well known; see, e.g.,  Lemmas~3 and 4 in \cite{Erdelyi}. 
 Its simple straightforward proof   can be found in \cite{II}.

\begin{corollary}\label{lma.d4}
Let the function $\sigma_*$ be given by \eqref{d6}, and let $h_{*}= {\sf L} \sigma_*$ be 
 its Laplace transform. Then
$$
h_{*}=\kappa_0 h_0+ b_\infty \kappa_\infty+\wt {h_{*}},
$$
where the model kernels $h_0$, $h_\infty$ are defined by 
\begin{equation}
h_0(t)=t^{-1}|\log t|^{-\alpha}\chi_0(t),
\quad
h_\infty(t)=t^{-1} |\log t |^{-\alpha}\chi_\infty(t)
\label{a7b}
\end{equation}
 and the error term $\wt {h}\in C^\infty(\Bbb R_+)$ satisfies the estimates
\[
|\wt {h_{*}}^{(m)}(t)|\leq C_m t^{-1-m}\la\log t\ra^{-\alpha-1}, \quad t>0,
\]
for all integers $m\geq0$.
\end{corollary}

 Our model operator is the Hankel operator $H_* = H (h_{*})$ with kernel $h_*= {\sf L} \sigma_*$. 

\begin{lemma}\label{lma.c0}
The eigenvalues of the  operator $H _*$
obey the asymptotic relation
\begin{equation}
\lambda_n^\pm(H_*)
=
a^\pm n^{-\alpha}+o(n^{-\alpha}), 
\quad 
n\to\infty,
\label{e11}
\end{equation}
where the coefficients $a^\pm$ are given by \eqref{eq:Al}. 
\end{lemma}

Indeed,
according to Theorem~\ref{SIG}, the Hankel operator $H_*$
is unitarily equivalent to the $\Psi$DO $A_{*} = v(X) s_*(D)v(X)$
in $L^2(\Bbb R)$. As usual,  $v (x)$ is the standard function \eqref{eq:vv}, and it follows from \eqref{eq:sig2} and \eqref{d6} that
$$
s_{*}(\xi)
=
\sigma_*(e^{-\xi})
=
\kappa_\infty |\xi|^{-\alpha}\chi_0(e^{-\xi})
+
\kappa_0 |\xi |^{-\alpha}\chi_\infty(e^{-\xi}), 
\quad \xi\in\Bbb R.
$$
This function   belongs to 
$C^\infty ({\Bbb R})$ and has the asymptotic behaviour \eqref{b2} with $s_{\infty}= \kappa_{\infty}$ 
and $s_{-\infty}= \kappa_0$. 
Therefore Theorem~\ref{thm.b2} (Weyl spectral asymptotics of $\Psi$DO) applies to the operator $A_{*}$.
This yields the asymptotic formula 
\begin{equation}
\lambda_n^\pm(A_*)
=
a^\pm
n^{-\alpha}+o(n^{-\alpha}), 
\quad 
n\to\infty,
\label{e11a}\end{equation}
where 
\[
a^\pm
=
(2\pi)^{-\alpha}
\bigl((\kappa_0)_\pm^{1/\alpha}+(\kappa_\infty)_\pm^{1/\alpha}\bigr)^\alpha
\biggl(\int_{-\infty}^\infty(\pi(\cosh(\pi x))^{-1})^{1/\alpha} dx\biggr)^\alpha.
\]
Using the change of variables $y=(\cosh(\pi x))^2$,     it is easy to check that the coefficients $a^\pm$ here and in \eqref{eq:Al} coincide. In view of Theorem~\ref{thm.e2} the operators $H_*$ and $A_*$ are unitarily equivalent so that
relation \eqref{e11a} yields  \eqref{e11}.
$\Box$

Now we are in a position to conclude the proof of Theorem~\ref{thm.d1}.
By its hypotheses, we have the representation 
$$
h=\kappa_0 h_0+ \kappa_\infty h_\infty+\wt {h },
$$
where $h_0$ and $ h_\infty$ are  given by \eqref{a7b} and  $\wt {h }$ satisfies the assumptions of Theorem~\ref{thm.a2} (singular value estimates). Therefore  it follows from
 Corollary~\ref{lma.d4}  that the difference 
$$
h- h_*= \wt {h }- \wt {h_{*} }
$$
also satisfies the hypothesis of Theorem~\ref{thm.a2} and hence
\begin{equation}
s_n(H -H_*) =o(n^{-\alpha}), \quad n\to\infty.
\label{eq:We}
\end{equation}
In view of the abstract Lemma~\ref{lma.b1} the asymptotic formula \eqref{d3a} is a direct consequence of \eqref{e11}
and   \eqref{eq:We}. 
\qed

   \medskip

 {\bf 4.4.}
 We emphasize that  the asymptotics of the spectrum of integral Hankel operators \e{eq:H1} is determined   by the behavior of $h(t) $  as $t\to 0$ and $t\to\infty$ as well as by local singularities  of $h(t)$. Following \cite{Y},
  let us consider Hankel operators whose integral kernels (or their derivatives) have   jumps of continuity at some positive point.
 
  \begin{theorem}\label{power}
Let $l\in{\Bbb Z}_{+}$, $t_{0}>0$, and let  $h(t)= h_{0} (t_{0}-t)^l $
 for $t\leq t_{0}$ and 
$h(t)=  0$ for $t > t_{0}$. Then eigenvalues of the Hankel operator $H$ have the asymptotics
\begin{equation}
\lambda_{n}^{\pm}(H)=   | h_{0}|  l!  (2\pi)^{-l-1}  t_{0}^{l+1} n^{-l-1} (1+ O(n^{-1}))
\label{eq:remm}
\end{equation}
as $n\to \infty$.
 \end{theorem}
 
 Of course, the exact expression for $h(t)$ is inessential. Indeed,   if a real function
 $v(t)$ satisfies the assumptions of   Theorem~\ref{thm.a2} with $\alpha=l$, then singular numbers $s_{n} (V)$
     of the Hankel operator  $V$ with kernel $v(t)$ satisfy the bound   $  s_{n}(V)= o(n^{- l -1})$. Therefore, in view of Lemma~\ref{lma.b1} asymptotics \e{eq:remm} remains true for the eigenvalues of the Hankel operator $H+V$; however,  in this case the remainder $O(n^{-1})$ in \e{eq:remm} should be replaced by $o(1)$.   

 We emphasize that,  according to \e{eq:remm}, the leading terms of the asymptotics of positive and negative eigenvalues of the Hankel operator  $H$ are the same.     Of course, if $h(t)$ becomes smoother ($l$ increases), then   eigenvalues   of $H$ decrease faster as $n\to \infty$.   Observe  that for $l=0$ (when the kernel itself is discontinuous), the Hankel operator $H$ does not belong to the trace class.
   We finally note that under the assumptions of Theorem~\ref{power} the asymptotics of singular values of the operator $H$ was found long ago in \cite{Part}.

\section{Discrete case} 

{\bf 5.1.}
 In the discrete case, the role of derivatives $h^{(m)}(t)$ of a function $h(t)$ is played by iterated differences $g^{(m)}(j)$ of a sequence  $g(j)$. Those are defined iteratively by setting $g^{(0)}(j)= g (j)$ and
 \[
 g^{(m)}(j)=g^{(m-1)}(j+1)-g^{(m-1)}(j),\q j\geq 0.
 \]
 
 The following result plays the role of Theorem~\ref{thm.a2}.

\begin{theorem}\label{thm.a1}\cite{I}
Let $\alpha>0$ and  let   $g$ be a sequence of complex numbers that satisfies
$$
g^{(m)}(j)=o(j^{-1-m}\la\log j\ra^{-\alpha}),\q j\to\infty,
$$
for all   $m=0,1,\ldots, N(\alpha)$ with $N(\alpha)$ defined by \e{eq:Mm}.
Then $s_n(G)=o(n^{-\alpha})$ as $ n\to\infty$.
\end{theorem}

Theorem~\ref{thm.a1} remains true if $o (\cdot)$ (on both occasions above)  is replaced by $O  (\cdot)$.
 
 Below is our main result in the discrete case.
\begin{theorem}\label{cr.a3}
Let $\alpha>0$, $\kappa_1, \kappa_{-1} \in\Bbb R$, 
and let $g(j)$ be a sequence of real numbers given $($for $ j\geq2)$ by 
\begin{equation}
g(j) =(\kappa_1+ (-1)^j \kappa_{-1} )j^{-1}(\log j)^{-\alpha}+\wt{g}_1(j)+(-1)^j \wt{g}_{-1}(j) 
\label{a2a}
\end{equation}
where the error terms $\wt{g}_{\pm 1}$ satisfy  the conditions of Theorem~\ref{thm.a1}. 
Then the eigenvalues of  the corresponding Hankel operator $G$ 
have the asymptotic behaviour
\[
\lambda_n^\pm(G)
=
b^\pm n^{-\alpha}+o(n^{-\alpha}), \q n\to\infty,
\]
where
\begin{equation}
b^\pm=\tau (\alpha)\bigl((\kappa_1)^{1/\alpha}_\pm+(\kappa_{-1})^{1/\alpha}_\pm\bigr)^\alpha
\label{e3}
\end{equation}
and $\tau(\alpha)$ is given by \eqref{a4}. 
\end{theorem}

\medskip

{\bf 5.2.} 
Let  us describe the plan of the proof of Theorem~\ref{cr.a3}. We follow the same steps as in Section~4, 
but instead of the Laplace transform $h_*={\sf L} \sigma_*$ of the function $\sigma(\lambda)$, $\lambda>0$, we consider the sequence of moments 
\begin{equation}
g_*(j)=\int_{-1}^1 \eta_*(\mu)\mu^j d\mu, \quad j\geq 0 ,
\label{e19}
\end{equation}
of  the function:
\begin{equation}
\eta_*(\mu)
=
\Big|\log \frac{1+\mu}{2(1-\mu)}\Big|^{-\alpha}
\left(
\kappa_1\chi_\infty\Big(\frac{1+\mu}{2(1-\mu)}\Big)
+
\kappa_{-1}\chi_0\Big(2\frac{1+\mu}{1-\mu}\Big)
\right) 
\label{e8}
\end{equation}
where the smooth cut-off functions $\chi_\infty$ and $\chi_0$ are given by equalities \eqref{a7a}.
Note that the function $\eta_*$ belongs to the class $C^\infty (-1,1)$ and
has the following asymptotic behaviour:
\begin{align*}
\eta_*(\mu)&=
\kappa_1 |\log(1-\mu)|^{-\alpha}+o(|\log(1-\mu)|^{-\alpha}),\quad \mu\to1,
\\
\eta_*(\mu)&=
\kappa_{-1} |\log(1+\mu)|^{-\alpha}+o(|\log(1+\mu)|^{-\alpha}),\quad \mu\to-1.
\end{align*}

These relations allow us to
 obtain the asymptotics of the sequence $g(j)$ as $j\to\infty$. 
We use again Lemma~\ref{L} but one needs
 to replace the continuous parameter $t$ with the discrete one $j$.
  The following assertion plays the  role of   Corollary~\ref{lma.d4}.

\begin{lemma}\label{lma.e3}
The sequence $g_*(j)$  defined by \eqref{e19}, \eqref{e8} has the  asymptotics 
\begin{equation}
g_{*}(j) =(\kappa_1+  (-1)^j \kappa_{-1})j^{-1}(\log j)^{-\alpha}+ \wt g_{1}(j)+(-1)^j \wt g_{-1}(j) 
\label{eq:a2a}
\end{equation}
where the error terms $\wt g_{\pm 1}(j)$  satisfy the estimates
\[
 \wt g_{\pm 1}^{(m)} (j)
=O(j^{-1-m}(\log j)^{-\alpha-1}), \quad j\to\infty,
\]
for all $m=0,1,2,\dots$. 
\end{lemma}

Let $v(x)$ be the   function \eqref{eq:vv} and
\[
s_{*}(\xi)=\eta_*\left(\frac{2e^{-\xi}-1}{2e^{-\xi}+1}\right), 
\quad \xi\in\Bbb R.
\]
Theorem~\ref{thm.e2} implies that the Hankel operator $G_*$ with matrix elements \eqref{e19}  is unitarily equivalent to the $\Psi$DO $A_{*} = v(X) s_*(D)v (X)$ acting
in $L^2(\Bbb R)$ .

By the definition \eqref{e8} of $\eta_* (\mu)$, we have
$$
s_{*}(\xi)=|\xi|^{-\alpha}(\kappa_1\chi_\infty(e^{-\xi})+ \kappa_{-1}\chi_0(4e^{-\xi})), 
\quad 
\xi\in\Bbb R.
$$
Applying Theorem~\ref{thm.b2} to the $\Psi$DO  $A_{*}$, we see that
\begin{equation}
\lambda_n^\pm(G_*)
=
\lambda_n^\pm(A_*)
=
b^\pm n^{-\alpha}
+
o(n^{-\alpha}),
\quad
n\to\infty,
\label{e18}\end{equation}
where 
 the numbers $b^\pm$ are   given by  \eqref{e3}.

\medskip

{\bf 5.3.}
Now we can conclude the
proof of Theorem~\ref{cr.a3}.
Comparing \eqref{a2a} and \eqref{eq:a2a}, we see that
\[
g(j) - g_{*} (j) = f_1(j)+(-1)^j f_{-1}(j) 
\]
where the error terms $f_{\pm 1}(j)= g_{\pm 1}(j)- \widetilde{g}_{\pm 1}(j)$ satisfy  the condition
$$
f_{\pm 1}^{(m)}(j)= o(j^{-1-m}  (\log j)^{-\alpha}), \quad j\to\infty,
$$
for all $m=0,1,\dots, N(\alpha)$. According to  Theorem~\ref{thm.a1} (singular value estimates) we have
$s_n(G(f_{\pm 1}))=o(n^{-\alpha})$ as $n\to\infty$.

Put $ \wt f_{-1}(j)=(-1)^j f_{-1}(j) $. Then   $s_n(G (\wt f_{-1}))=s_n(G(f_{-1}))$ and hence
\begin{equation}
s_n(G(f_1 + \wt f_{-1}))=o(n^{-\alpha}), \quad n\to\infty.
\label{e14+}
\end{equation}
In view of  \eqref{e18}  and \eqref{e14+}, we can apply the abstract Lemma~\ref{lma.b1} to the operators $K_{0}= G(g_*)$ and $K=G (f_1 + \wt f_{-1})$ . This
 yields the eigenvalue asymptotics \eqref{e3} for the operator $G=G(g)=K_{0}+K$. 
\qed

  \section{Generalizations and applications}
  
 {\bf 6.1.} 
 In this subsection we state results on asymptotic behavior of singular values of Hankel operators. Now the operators are not assumed to be self-adjoint. Recall that the numerical coefficient 
$\tau(\alpha)$ is given by formula \eqref{a4}.

 It is convenient to start with the discrete case.
 
\begin{theorem}\label{thm.a4}\cite[Theorem 3.1]{III}
Let $\alpha>0$, let $\zeta_1, \dots,\zeta_L\in\Bbb T$ be pairwise distinct numbers, 
and let $\kappa_{1},\dots, \kappa_L\in\Bbb C$.
Let $g(j)$ be a sequence of complex numbers such that
\begin{equation}
g(j)=\sum_{\ell=1}^L \bigl(\kappa_\ell j^{-1}(\log j)^{-\alpha}+ \wt{g}_\ell(j)\bigr)\zeta_\ell^{-j},
\quad
j\geq2,
\label{a13}
\end{equation}
where the error terms $\wt{g}_\ell$, $\ell=1,\ldots, L$, satisfy the estimates
\begin{equation}
\wt{g}_\ell^{(m)}(j)
=o(j^{-1-m}(\log j)^{-\alpha}), 
\quad j\to\infty, 
\label{a5a}
\end{equation}
for all   $m=0,1,\dots, N(\alpha)$
with $ N(\alpha)$  given by \eqref{eq:Mm}.
Then the singular values of the Hankel
  operator $G $ defined in $\ell^2({\Bbb Z}_+)$   by formula \eqref{eq:Hd}  satisfy the asymptotic relation
 \begin{equation}
s_n(G)= b \, n^{-\alpha}+o(n^{-\alpha}), 
\quad n\to\infty,
\label{eq:Lm}
\end{equation}
where
\begin{equation}
b =\tau(\alpha)\Big(\sum_{\ell=1}^L |\kappa_\ell |^{1/\alpha}\Big)^\alpha .
\label{eq:Lm1}
\end{equation}
\end{theorem}

The plan of the proof of Theorem~\ref{thm.a4}  is the following. For $L=1$,   Theorem~\ref{thm.a4} is a consequence of Theorem~\ref{cr.a3} for the particular case $\kappa_{-1}=0$. To pass to the general case, one needs the notion of the symbol $\omega(\mu)$ of a Hankel
  operator $G$. The function $\omega(\mu)$ can be defined by the relation 
   \begin{equation}
g(j)= \int_{\Bbb T} \omega(\mu) \mu^{-n} d{\bf m}(\mu)
\label{eq:ssy}\end{equation}
where $d{\bf m}(\mu)$ is the normalized Lebesgue measure on the unit circle ${\Bbb T}$. Of course, the function  $\omega(\mu)$ satisfying \e{eq:ssy} is not  unique.

Consider  the leading part  
\begin{equation}
g_{\,\rm lead} (j)=\sum_{\ell=1}^L  \kappa_\ell j^{-1}(\log j)^{-\alpha} \zeta_\ell^{-j},
\quad
j\geq2,
\label{eq:a13X*}
\end{equation}
of sequence \e{a13}.  It can be easily checked that for every $\alpha>0$ the function
   \[
  \omega_{0}(\mu) = \sum_{j=2}^\infty j^{-1}(\log j)^{-\alpha}(\mu^j-\bar{\mu}^j), \q \mu\in {\Bbb T} ,
\]
is bounded and $\omega_{0}\in C^\infty ( {\Bbb T}\setminus\{1\}$). This means that its  singular support is
$\sing\supp \omega_{0} = \{1\}$. 
Therefore the singular support of  the symbol $\omega_{\,\rm lead} (\mu)$ corresponding to $g_{\,\rm lead} (j)$ consists of the points  $\z_{1}, \ldots, \z_{L}$.
It can be deduced from this property that the singular value counting function
$
\# \{n: s_n (G_{\rm lead})>\varepsilon\}$ 
 of the Hankel operator $G_{\,\rm lead}$  with the matrix elements $g_{\,\rm lead} (j)$ is asymptotically (as $\varepsilon\to + 0$)  the sum of such functions for separate terms in \e{eq:a13X*}. This fact is called the localization principle in \cite{III}. In terms of singular values the result on counting functions is equivalent to relations \e{eq:Lm}, \e{eq:Lm1} for the operator $G_{\,\rm lead}$. The singular values of the operator $G-G_{\,\rm lead}$ can be easily estimated with the help of Theorem~\ref{thm.a1}. $\Box$

 
In the continuous case, we have the following result.

\begin{theorem}\label{thm.d3}\cite[Theorem 5.1]{III}
Let $\alpha>0$, let $\rho_1,\dots, \rho_L\in\Bbb R$ 
be pairwise distinct numbers and let $\kappa_0, \kappa_1,\dots, \kappa_L\in\Bbb C$. 
Let the number $N(\alpha)$ be given by \eqref{eq:Mm}. 
Suppose that $h\in L^\infty_{\rm loc}({\Bbb R}_+)$ if $\alpha<1/2$ and $h\in C^{N(\alpha)} ({\Bbb R}_+)$ if $\alpha\geq1/2$. 
Assume that 
\begin{align*}
h(t)&=\sum_{\ell=1}^L \bigl( \kappa_\ell t^{-1}(\log t)^{-\alpha}+  {\wt h}_\ell(t)\bigr)e^{- i \rho_\ell t}, 
\quad t\geq2,
\\
h(t)&=\kappa_0t^{-1}\bigl(\log(1/t)\bigr)^{-\alpha}+  {\wt h}_0(t), 
\quad t\leq 1/2,
\end{align*}
where the error terms ${\wt h}_\ell$ and their derivatives ${\wt h}_\ell^{(m)}$
satisfy the estimates
   \begin{equation}
{\wt h}_\ell^{(m)}(t)
=
o( t^{-1-m}\la\log t\ra^{-\alpha}), 
\quad 
m=0,\dots,N(\alpha),  
\label{z4}
\end{equation}
as $ t\to\infty$ for $\ell=1,\ldots, L$ and
 as $t\to 0$ for $\ell=0$.
Then the singular values of the integral Hankel
operator $H$ with kernel $h(t)$ in $L^2({\Bbb R}_+)$  satisfy the asymptotic relation
\[
s_n(H)= a \, n^{-\alpha}+o(n^{-\alpha}), \quad n\to\infty,
\]
where 
\[
a  =\tau(\alpha)\Big(\sum_{\ell=0}^L | \kappa_\ell|^{1/\alpha}\Big)^\alpha.
\]
\end{theorem}

The proof of Theorem~\ref{thm.d3} follows the same general outline
as that of Theorem~\ref{thm.a4}. In the continuous case the symbol may be defined by the 
  relation 
   \[
h(t)= \int_{-\infty}^\infty \Omega(x) e^{-ixt} dx.
\]
Compared with the proof of Theorem~\ref{thm.d3} the only essential difference is
  that the singularity of the kernel $h(t)$ at 
$t=0$ has to be treated separately. It corresponds to the singularity
of the symbol $\Omega(x)$ at infinity.

\medskip

{\bf 6.2.}
Here we   find
asymptotics of   eigenvalues  of self-adjoint Hankel operators. The first result generalizes Theorem~\ref{cr.a3}.
Now we  consider the real sequences of the form \eqref{a13}.



\begin{theorem}\label{thm.a5}\cite[Theorem 5.7]{V}
Let $\alpha>0$, $p=1/\alpha$; let $\zeta_1,\dots,\zeta_L\in\Bbb T$ be pairwise distinct points
with $\Im \zeta_\ell>0$, and let $\varkappa_{1}, \varkappa_{-1}\in\Bbb R$, 
$\kappa_1,\dots, \kappa_L\in\Bbb C$. 
Let $g (j)$ be a sequence of real numbers such that 
\begin{multline}
g(j)=
\varkappa_1j^{-1}(\log j)^{-\alpha}
+ \wt{\sf g}_{1}(j) + (-1)^{j} \big(
\varkappa_{-1} j^{-1} (\log j)^{-\alpha} +  \wt{\sf g}_{-1}(j) \big)
\\
+
2\Re \sum_{\ell=1}^L   \zeta_\ell^{-j} \big( \kappa_\ell  j^{-1}(\log j)^{-\alpha}
+ \wt{g}_{\ell}(j) \big),
\quad 
j\geq2,
\label{a17}
\end{multline}
where all error terms $\wt{\sf g}_{1}, \wt{\sf g}_{-1},\wt{g}_{1}, \dots, \wt{g}_{L}$   
obey condition  \eqref{a5a} for  $m=0,1,\dots,N (\alpha)$ $(N(\alpha)$ is given by \eqref{eq:Mm}$)$. 
Then the eigenvalues of the Hankel operator $G$ with matrix elements \e{a17}
satisfy the asymptotic relation \e{e3}
with the coefficient $b^\pm$  defined by
\[
b^\pm
=
\tau (\alpha)\bigl((\varkappa_{-1})_\pm^{1/\alpha}+(\varkappa_{1})_\pm^{1/\alpha} 
+\sum_{\ell=1}^L | \kappa_\ell|^{1/\alpha}\bigr)^\alpha.
\]   
\end{theorem}

Compared to the proof of Theorem~\ref{thm.a4} one has to additionally use the so called symmetry principle (see \cite{V}). It states that if the singular support of the symbol of a compact self-adjoint Hankel operator $G$  does not contain the points $1$ and $-1$, then the spectrum
of $G$ is asymptotically symmetric with respect to the point zero.

In the continuous case, we consider real kernels $h(t)$ that are singular at $t=0$ and contain several oscillating terms at  infinity. The assertion below plays the role of Theorem~\ref{thm.a5}, and its proof follows essentially the same lines.
Similarly to the proof of  Theorem~\ref{thm.d3}, the contributions of the points $t=\infty$ and $t=0$ should be considered separately.

\begin{theorem}\label{thm.a7}\cite[Theorem 6.6]{V}
Let $\alpha>0$, let $\rho_1,\dots,\rho_L  \in{\Bbb R}_{+}$ be pairwise distinct   numbers, 
and let $\kappa_0, \kappa_\infty\in {\Bbb R}$, $\kappa_1,\dots, \kappa_L\in {\Bbb C}$. 
Suppose that $h\in L^\infty_{\rm loc} ({\Bbb R}_{+})$ if $\alpha<1/2$ and $h\in C^{N(\alpha)} ({\Bbb R}_{+})$  if $\alpha\geq 1/2$  where the number $N(\alpha)$ is   given by \eqref{eq:Mm}.
Assume  that
\[
h(t)= \kappa_\infty  t^{-1}(\log t)^{-\alpha}+\wt{h}_\infty(t)
+2\Re \sum_{\ell=1}^L 
\bigl(\kappa_\ell  t^{-1}(\log t)^{-\alpha}+\wt{h}_\ell(t)\bigr)e^{- i \rho_\ell t}
\]
  for $ t\geq2 $ and
\[
h(t)= \kappa_0t^{-1}\bigl(\log(1/t)\bigr)^{-\alpha}+  \wt{h}_0(t), 
\quad t\leq 1/2,
\]
where the error terms $\wt{h}_\infty$, $\wt{h}_1, \ldots, \wt{h}_L$ obey the estimates
\eqref{z4}
as $ t\to\infty$  and $\wt{h}_0$ obeys these estimates
 as $t\to 0$. 
Then the eigenvalues of the integral Hankel operator $H$ with kernel  $h(t)$
satisfy  asymptotic relation 
\e{d3a}
where 
\[
a^\pm 
=
\tau(\alpha)\bigl((\kappa_0)_\pm^{1/\alpha}+(\kappa_\infty )_\pm^{1/\alpha} +\sum_{\ell= 1}^L |\kappa_\ell|^{1/\alpha}\bigr)^\alpha.
\]
\end{theorem}

\medskip

{\bf 6.3.}
As an application of Theorem~\ref{thm.a4},
  we  state here  a result on rational approximations
    in the $\BMO$-norm  of    functions  $\omega(z)$ (of bounded mean variation on the unit circle $\Bbb T$) analytic on the unit disc $\Bbb D$ and acquiring some singularities on  $\Bbb T$. We are interested in singularities of logarithmic type.
   We study the asymptotic behavior as $n\to\infty$ of the distance $\dist_{\BMO} \{\omega, {\cal R}_{n}\}$ in the $\BMO$-norm between $\omega$
    and  the set  ${\cal R}_{n}$ of all rational functions of degree $\leq n$ without poles 
on $\ov{\Bbb D}$.
       A short description of relevant results in this vast domain can be found in \cite{IV}.
 In view of the Adamyan-Arov-Kre\u{\i}n theorem the problem considered is equivalent to the
study of the asymptotic behaviour of singular values of the Hankel operator with symbol $\omega (z)$.

Let us describe the class of admissible functions $\omega(z)$.
Let $u(z)$ be analytic in $\Bbb D$, $u\in C^\infty(\ov{\Bbb D})$; fix some $\zeta = e^{i\varphi }\in\Bbb T$ and
assume that 
\begin{equation}
-\log ( \zeta - z)+u (z)\neq 0, \quad z\in \ov{\Bbb D}.
\label{LA2a}
\end{equation}
Define
\[
\omega (z)=\bigl(-\log ( \zeta - z)+u(z)\bigr)^{1-\alpha},\quad z \in\Bbb D, \quad \alpha>0.
\]
We have introduced $u(z)$ to avoid irrelevant singularities of $\omega (z)$ inside $\Bbb D$.
The branch of the analytic function $\log ( \zeta - z)$ is fixed by the condition 
$\log ( \zeta - z)=\log ( 1-r)+i\varphi$ if $z=r \zeta$ for $r\in (0,1)$. 
We fix $\arg \bigl(-\log ( \zeta - z)+ u (z) \bigr)$ by the condition that it tends to zero as $z=r e^{i\varphi}$ and $r\to 1-0$. Obviously, the function $\omega(z)$ is analytic in the  unit disc $\Bbb D$ and is smooth up to the boundary $\Bbb T$, except at the point $z=\zeta$. 

Theorem~\ref{thm.a4} allows us to consider $\omega(z)$  as well  as finite sums of such functions.

\begin{theorem}\label{Analytic}\cite[Theorem 3.8]{IV}
Let   $\zeta_1, \zeta_2,\dots,\zeta_L\in\Bbb T$ be pairwise distinct points, and  let functions   $v_{\ell}, u_{\ell}$, $\ell=1,\ldots, L$, be analytic in     $\Bbb D$ and  $v_{\ell}, u_{\ell}\in C^\infty(\ov{\Bbb D})$. Assume that \eqref{LA2a} is satisfied for all $u_\ell$, $\zeta_\ell$ and set
$$
\omega (z)=\sum_{\ell=1}^L v_{\ell} (z)\bigl(-\log ( \zeta_{\ell}- z)+u_{\ell} (z)\bigr)^{1-\alpha},  \quad \alpha>0 .
$$
Then there exists the limit
$$
\lim_{n\to\infty}n^\alpha
\dist_{\BMO} \{\omega, {\cal R}_{n}\}
= |1-\alpha| \tau(\alpha)\Bigl(\sum_{\ell=1}^L |v_\ell (\zeta_\ell)| ^{1/\alpha}\Bigr)^\alpha .
$$
\end{theorem} 

Note  that for $\alpha<1$, the functions $\omega (\z)$ are unbounded as $\z\in{\Bbb T}$ tends to one of the points $\zeta_{\ell}$ so that their approximation in the norm of $C({\Bbb T})$ is  a priori impossible.


 \end{document}